\documentclass[11pt]{article}
\usepackage{amsmath, amscd, amsfonts, amssymb, amsthm, tikz, times}
\newtheorem{thm}{Theorem}
\newtheorem{prop}{Proposition}

\newtheorem{cor}{Corollary}

\newcommand{\lar}{\longrightarrow}
\newcommand{\ff}{\infty}
\newcommand {\rrr}[1]{(\ref{#1})}
\newcommand{\RR}{{\mathbb R}}
\newcommand{\vski}{\vspace{11pt}}

\title{Sum rules for effective resistances in infinite graphs}
\author{Greg Markowsky\\School of Mathematical Sciences\\Monash University, Melbourne, Australia\\and\\Jos\'e Luis Palacios\\ Department of Electrical and Computer Engineering\\  The University of New Mexico, Albuquerque, USA.}
\begin{document}

\maketitle
\vskip .5 in
\begin{abstract}
Extending work of Foster, Doyle, and others, we show how the Foster Theorems, a family of results concerning effective resistances on finite graphs, can in certain cases be extended to infinite graphs. A family of sum rules is then obtained, which allows one to easily calculate the sum of the resistances over all paths of a given length. The results are illustrated with some of the most common grids in the plane, including the square, triangular, and hexagonal grids.
\end{abstract}

\vskip .2 in
Running title: Effective resistances via Foster's formulas
\vfill
\eject
\section{Introduction}

Evaluating the electric resistance between points in infinite grids of resistors is a classical problem (see \cite{van}) which has retained its interest for modern researchers (for instance \cite{A}, \cite{CSD}). For example, the following resistances from a reference point $0$ have been calculated on,
respectively, the square, triangular, hexagonal, and truncated square tilings (\cite[Ch. XVI.7]{van},\cite{A}, \cite{CSD}).

\vspace{-1cm}
\begin{center}
\includegraphics[height=5in]{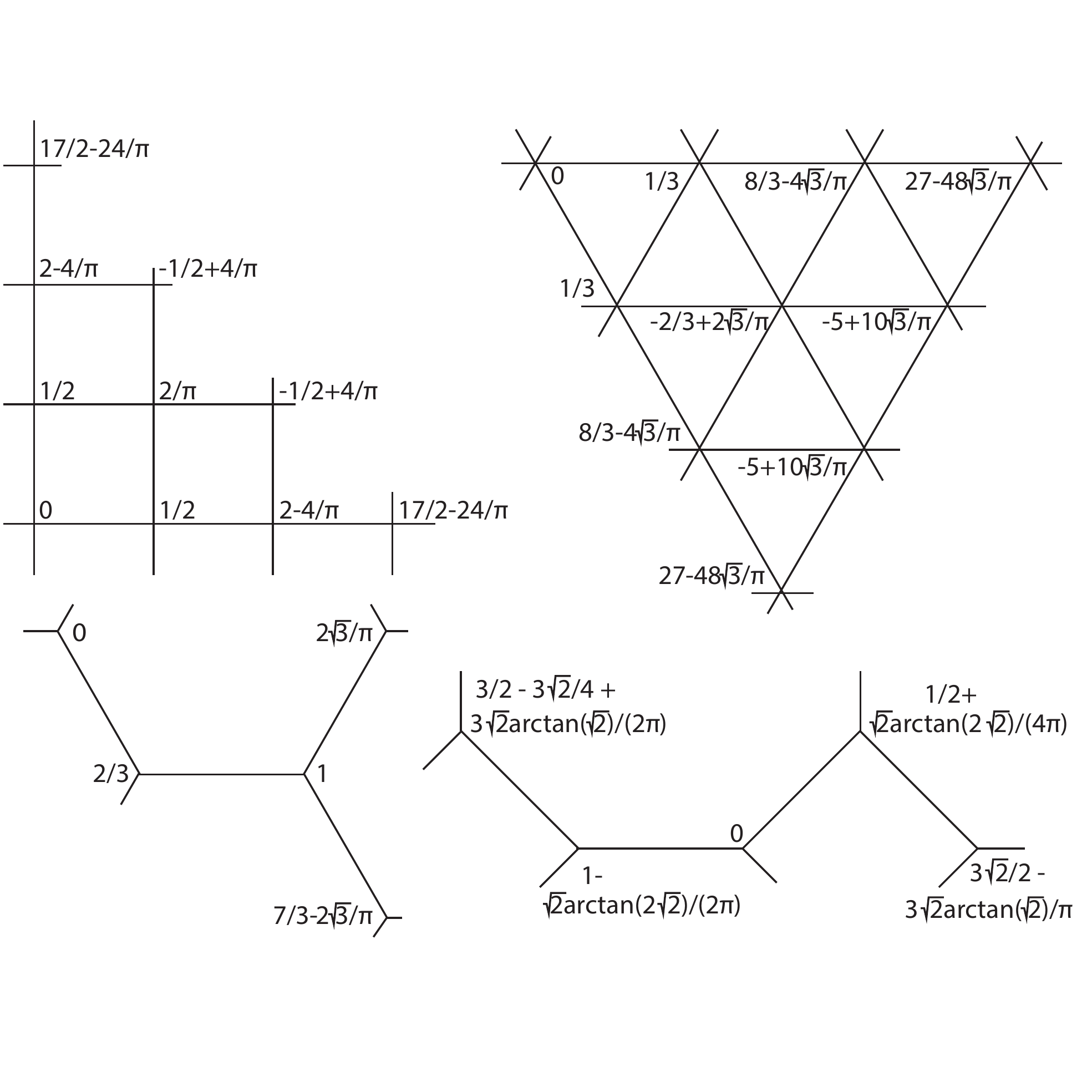}
\end{center}
Note that, taking into account the symmetries of the lattices, the first three pictures show the resistances between all pairs of vertices on the respective grids of distance at most 3, while the last picture covers all pairs
of points of distance at most 2.
Several interesting features may be apparent in examining these quantities. To begin with, we note that a number of the resistances are simple, rational numbers, specifically the resistances in points adjacent to $0$ in the first
three lattices, as well as the point of distance $2$ from $0$ in the hexagonal lattice. On the other hand, all others given are irrational, and are in fact not even algebraic. This may initially appear to be no more than a curiosity,
however a closer examination yields an interesting fact, namely that there are various linear combinations of effective resistances, with positive integer coefficients, which yield positive integers. For instance, in the square grid, if we consider only
the points of distance $2$ from $0$, we see that twice the resistance to $(1,1)$ plus the resistance to $(2,0)$ equals 2. Similarly, the sum of the resistances of the two different types of points of distance 3
from 0 in the hexagonal grid is rational.
The main motivation of this paper is to explain this phenomenon in terms of the Foster Theorems, which is a series of formulas concerning the average resistances between points of a given distance in a finite graph. We will show
how, in many cases, Foster's Theorems can be extended to infinite graphs, yielding simple summation rules for resistances as described in the title.

\vski

In order for this be realized, the infinite graph must (a) possess a great deal
of symmetry, and (b) satisfy a condition on the growth of boundaries of large balls. We begin with a discussion of symmetry, postponing the growth condition until Section \ref{pink}. In what follows we will employ the notation and terminology of mathematicians, but will also attempt to translate the notions into plain English in the hope of being accessible to a wider audience. A {\it graph} is a set of vertices and edges between them. If $x,y$ are vertices with an edge between them, then we will say that $x$ is {\it adjacent} to $y$, and write $x \sim y$. An {\it automorphism} of a graph is a function $\phi$ which maps the vertex set of the graph to itself such that (a) $\phi$ is one-to-one ($\phi(x)=\phi(y)$ only when $x=y$), (b) $\phi$ is onto (for any $y$, there is an $x$ such that $\phi(x)=y$), and (c) if $x \sim y$ then $\phi(x) \sim \phi(y)$. For example, automorphisms of the planar grids shown in the image above are simply Euclidean transformations (translations, rotations, and reflections) under which the lattices are invariant. A graph is {\it vertex transitive} if, given any two vertices $x,y$ there is a graph automorphism taking $x$ to $y$, and {\it s-arc
transitive} if, given any two paths $a,b$ of length $s$ (which do not repeat vertices) there is a graph automorphism taking $a$ to $b$; in the case $s=1$ we will say the graph is {\it edge transitive}. Put more intuitively, a graph is vertex transitive if all of the vertices "look the same", edge transitive if the same statement holds for all edges, and s-arc transitive if the same holds for all paths of length $s$. These types of transitivity need not go together, although not infrequently they do. For example, all of the graphs in the image above are vertex
transitive,  and all but the last are edge transitive. The hexagonal grid alone is 2-arc transitive, since every path of length 2 lies on a unique hexagon; note that it is not 3-arc transitive, though, since a path of length 3 may or may not lie on one hexagon. Graphs which are edge-transitive but not vertex transitive also exist; perhaps the easiest example is to place a vertex on the midpoint of each edge in the trangular, square, or hexagonal lattice (this is discussed more in Section \ref{split}).

\vski

All of the graphs we will consider will possess one of these types of symmetry. The paper will be organized as follows. In the next two sections, we will discuss Foster's Theorems on finite graphs, and work done by other authors on extending the theorems to infinite graphs. Sections \ref{vertran} and \ref{split} contain the
new results, and show how the theorems can be extended further to certain vertex transitive and $s$-arc transitive infinite graphs, yielding the sum rules referred to in the title of this paper. The proofs of the results are placed in an appendix at the end, which contains a few final remarks as well.

\section{Foster's Theorems}

Let $G$ be a finite graph whose edges are endowed with unit resistances, and let $d(x,y)$ denote the distance of the shortest path in the graph between vertices $x$ and $y$.
Let $n=|G|$ be the number of vertices of $G$, and $m$ be the number of edges of $G$. For any vertex $v$, let $deg(v)$ denote the degree of $v$, which is the number of vertices adjacent to $v$. Let $R_{xy}$ denote the electric resistance between the points $x,y$, whose standard definition can
be found in many sources, including \cite{DS}. R. M. Foster proved in \cite{F1} and \cite{F2} the two well-known formulas for finite graphs:
\begin{equation}
\label{foster1}
\sum_{v \in G} \sum_{y \sim v} R_{vy}=2(n-1)
\end{equation}
and
\begin{equation}
\label{foster2}
\sum_{v \in G} \sum_{x \sim y \sim v} {{R_{xv}}\over {deg(v)}}=2(n-2).
\end{equation}
It should be noted that these formulas are different than the ones normally given, because as we have written it the resistance between each pair of points is counted twice, introducing a factor of 2 on the right side; the reason
for presenting the formulas
in this way is that they extend more naturally to the infinite setting, as will be seen in later sections. In \cite{T2} a probabilistic interpretation and proof of Foster's First \rrr{foster1}
was given, and the ideas of that paper were extended in \cite{P3} to prove Foster's Second \rrr{foster2}. Foster's $r$-th Theorem was also alluded to in that work, and proved for $r=3$, before a complete statement was made and
proved in \cite{P4}; it should be noted that the result was discussed also in different contexts in \cite[Thm. G]{klein2002} and \cite[Prop. 2.3]{bendito2008}. In order to state the theorem, let us use $P$ to denote the $n \times n$ matrix with entries

\begin{equation} \label{}
P_{ij} = \left \{ \begin{array}{ll}
\frac{1}{deg(i)} & \qquad  \mbox{if } i \sim j  \\
0 & \qquad \mbox{if } i\nsim j\;,
\end{array} \right.
\end{equation}

In probabilistic terms, this is the transition matrix of the random walk on the graph (see the appendix for more on this aspect of the result). When written with an exponent, such as $P^s$, we mean $P$ raised to the $s$ power by matrix multiplication. We will denote the trace of a matrix $A$, which is the sum of the diagonal elements of $A$, as $tr(A)$. With this notation, Foster's $r$-th is as follows.

\begin{thm} \label{mpbaby}

\begin{equation} \label{marcopolo}
\sum_{i \in G} \sum_{j \sim v_{r-1} \sim \ldots \sim v_1 \sim i}\frac{R_{ij}}{deg(v_1)deg(v_2) \ldots deg(v_{r-1})} = 2\Big( \sum_{s=0}^{r-1} tr(P^s) - r\Big).
\end{equation}
\end{thm}

\noindent In fact, the result in \cite{P4} is more general than this, and given in a more complex form, but the form given here is sufficient for purposes; a few comments about the derivation of this result and the more general form can be found in the appendix. Note that $tr(P^0) = n$ and $tr(P^1)=0$, which show that \rrr{marcopolo} reduces to \rrr{foster1} and \rrr{foster2} in the cases $r=1,2$.

\section{Foster's Theorems and electric resistance on infinite graphs - a summary of the literature} \label{pink}

In what follows, we will at times refer to various infinite graphs embedded in Euclidean space as {\it grids}
or {\it lattices}, in accordance with common terminology found in the literature.
Foster used his formulas in \cite{F2} in the computation of the effective resistance between neighboring vertices of infinite grids arguing thus: if the graph is vertex and edge transitive, so that all effective resistances
between neighbors are the same, then (\ref{foster1}) can be rewritten as
\begin{equation}
\label{foster3}
R_{ij}={{2(n-1)}\over {kn}},
\end{equation}
where $k$ is the degree of the graph. As $n \rightarrow \infty$, this implies that in an infinite $k$-regular graph with edge symmetry we have

\begin{equation}
\label{foster4i}
R_{ij}={2\over k},
\end{equation}
whenever $d(i,j)=1$. If we assume instead the graph is 2-arc transitive then (\ref{foster2}) can be reduced to
\begin{equation}
\label{foster5}
R_{ij}={{2(n-2)}\over {(k-1)n}}.
\end{equation}
As $n\rightarrow \infty$, this implies that in an infinite 2-arc transitive graph we have
\begin{equation}
\label{foster6}
R_{ij}={2\over {k-1}},
\end{equation}
whenever $d(i,j)=2$.
Foster goes on to exemplify his results with the hexagonal lattice and the square lattice.  For the first the formulas (\ref{foster4i}) and (\ref{foster6}) imply that $R_{ij}={2\over 3}$ when $d(i,j)=1$ and $R_{ij}=1$ when
$d(i,j)=2$.  For the second, (\ref{foster4i}) implies  $R_{ij}={1\over 2}$ when $d(i,j)=1$, but (\ref{foster6}) cannot applied, since the square grid is not 2-arc transitive (some pairs of adjacent edges have endpoints contained in a single
square, while others do not).

A different and successful way to compute the effective resistance between points at an arbitrary distance in the square lattice appears in \cite{Ven} and \cite{van}, and consists of applying the superposition principle to
currents that are injected at a vertex and let to come out at infinity, or injected at infinity and let out at a vertex.  This method relies heavily on the symmetry of the grid and the mathematical problem involves the solution
of an infinite set of linear, inhomogeneous  difference equations which are solved by the method of separation of variables.  In \cite{A} this method of superposition of currents was improved and extended to other grids with high degree of symmetry: the triangular and hexagonal
lattices in two dimensions, and also to infinite cubic and hypercubic lattices in three and more dimensions. The results found by this method agree in these cases with the results obtained by Foster's Theorems. However, as elegant as the results and justifications involving the Foster Theorems are, the application of the theorems to infinite graphs must be made rigorous. To illustrate that there really is something worth worrying about here, note
that the reasoning given above, if applied to the infinite homogeneous binary
tree (which is the infinite graph with no cycles - that is, no non-self-intersecting paths with the same beginning and end point -  where every vertex has degree 3), would yield a resistance of $\frac{2}{3}$ between adjacent points; however, for any tree the resistance between adjacent points is $1$, since there is only one path between the points for the electricity to flow along. In fact, the reasoning can be made rigorous in many important
situations, as Doyle notes in \cite{Doy}, a survey which contains also many historical details.  In particular, Foster's arguments work for infinite graphs which have bounded degree and are {\it smallish}, meaning that they possess
a sequence of subgraphs $G_m$ such that $$(i)\: G_m \subseteq G_{m+1} \mbox{ and } \cup_m G_m = G, \mbox{ and }$$  $$(ii) \qquad \lim_{m \to \ff} {{|{\rm boundary}(G_m)|}\over {|G_m|}}=0.$$
Intuitively, this means that as $m$ gets large almost all of the points in $G_m$ should be away from the boundary of $G_m$. To illustrate, if we let $G$ be the triangular, square, or hexagonal lattice, and $G_m$ the intersection of $G$ with a disk of radius $m$ centered at the origin, then for large $m$ the number of points on the boundary of $G_m$ will behave like the circumference of the disk, while the number of points away from the boundary will behave like the area of the disk, so that the quotient in $(ii)$ will be something like $\frac{2 \pi m}{\pi m^2} \lar 0$, and we conclude that the grids are smallish. Following Doyle, we will call a sequence $G_m$ satisfying $(i)$ and $(ii)$ a {\it swelling sequence}.  Essentially, Doyle's argument boils down to the observation that for smallish graphs when $m$ becomes large most vertices and edges in $G_m$ are
far from the boundary. As such, the resistance across most edges in the subgraph
$G_m$ will be close simultaneously to the resistance across the same edge in the infinite grid as well as to the average resistance in $G_m$ given by Foster's Theorem. Letting $m \lar \infty$ completes the proof. Note that the
infinite homogeneous binary tree is not smallish, since the boundary of a large ball contains about half of the vertices in the ball. Doyle also departs from the extremely symmetric examples in the other works cited: he shows that
in an infinite grid that is edge transitive and smallish but not vertex transitive, for any neighboring vertices $i$ and $j$ one has

\begin{equation}
\label{doyle}
R_{ij}={{deg(i)+deg(j)}\over {deg(i)deg(j)}}.
\end{equation}
This clearly reduces to (\ref{foster4i}) if $deg(i)=deg(j)$, but also extends (\ref{foster4i}), since it is possible for edge transitive graphs to have two different types of vertices which have different degrees (for example, the
decorated hexagonal lattice discussed in Section \ref{split}).


\section{Vertex transitive grids} \label{vertran}

This section begins the new results, and we begin by showing that Foster's Theorems can say a great deal when a grid $G$ is vertex transitive. In this case, $G$ has a well defined degree $k$, and we will use the notation $\Delta_s$
to denote the number of paths of length $s$ which start and end at a given point; note that when calculating $\Delta_s$, one must count all paths of length $s$, even those which repeat vertices. We
also define $\Delta_0 = 1$. We then have

\begin{thm} \label{jess}
Let $G$ be a smallish, vertex-transitive, infinite graph of degree $k$. Then, for any vertex $v$ in $G$ and positive integer $r$, we have

\begin{equation} \label{biggie}
\sum_{j = v_r \sim v_{r-1} \sim \ldots \sim v_0 = v} R_{ij} = 2k^{r-1} \sum_{s=0}^{r-1} \frac{\Delta_s}{k^s}.
\end{equation}

\end{thm}

\noindent The proof can be found in the appendix (although the reader may like to compare directly with Theorem \ref{mpbaby}). As indicated above, $\Delta_0 = 1$; this is seen to be the correct definition by \rrr{zhang}, because $P^0$ is the identity matrix. Note also that $\Delta_1 = 0$. \rrr{biggie} can therefore be written as

\begin{equation} \label{biggie2}
\sum_{j = v_r \sim v_{r-1} \sim \ldots \sim v_0 = i} R_{ij} = 2k^{r-1} \Big(1 + \sum_{s=2}^{r-1} \frac{\Delta_s}{k^s}\Big),
\end{equation}
\noindent where the sum on the right side is taken to be 0 when $r \leq 2$. Let us isolate the first few cases of $r$ as corollaries. Foster's First Rule is as follows:

\begin{cor} \label{}
$$\sum_{x \sim v} R_{vx} = 2.$$
\end{cor}
\noindent Note that, if $G$ is edge transitive, then this implies that the resistance across any edge is $\frac{2}{k}$, which agrees with Foster's reasoning above. This applies and gives the correct values for the triangular, square, and hexagonal grids, but does not apply to the truncated square grid (bottom right in the picture at the beginning of the paper), which is not edge-transitive. In this case, we may still verify that $1-\sqrt{2} \arctan (2\sqrt{2})/2\pi + 2(1/2 + \sqrt{2} \arctan (2\sqrt{2})/4\pi) = 2$; note that the any vertex will have two edges attached to it which lie between a square and an octagon, which is the reason for doubling that resistance. Foster's Second is as follows:

\begin{cor} \label{}
$$\sum_{y \sim x \sim v} R_{yv} = 2k.$$
\end{cor}
\noindent Let us examine in more detail how this is applied. For the hexagonal lattice, say, there are two possibilities for a walk of length 2 starting at the origin: we may go out along an edge and then return to the origin, or we may go out along an edge and then go further from the origin in the second step. We can disregard the path that returns to the origin, since this returns a resistance of $0$ between the endpoints of the path, and the resistances between any pairs of points of distance 2 will be the same, since any path of length 2 lies on the edge of a unique hexagon and can therefore be mapped to each other by an automorphism. There are 6 such paths of length 2 emanating from any point: 3 choices for the first direction, then 2 for the second. Here $k$ is 3, so the corollary tells us that the sum of the resistances between the endpoints of these 6 paths is 6, or in other words that the resistance between any two points of distance 2 is 1. A more involved example is the triangular lattice. There are three types of paths of length 2 (not counting the ones where we return to the origin): we may move twice in the same direction, we may turn $60^\circ$ after the first step, or we may turn $120^\circ$ after the first step, which in fact brings us to a point adjacent to the origin. There are only 6 paths in which the two steps lie in the same direction, but 12 for each of the other types (since we may turn to either the right or left). Here $k$ is 6, and we obtain

\begin{equation} \label{}
12 \Big(\frac{1}{3}\Big) + 12 \Big(-\frac{2}{3} + \frac{2\sqrt{3}}{\pi}\Big) + 6 \Big(\frac{8}{3} - \frac{4\sqrt{3}}{\pi}\Big) = 12.
\end{equation}

\noindent

Similar calculations hold for the square and truncated square grids given in the picture. As is the case of Foster's Second for finite graphs, it is often easier to count paths of length 2 from the midpoint of the path rather than from an endpoint. Counting from an endpoint essentially counts each path twice, while
from the midpoint only counts it once. The following is therefore an equivalent statement for Foster's Second:

\begin{cor} \label{}
$$\sum_{x,y \sim v} R_{xy} = k.$$
\end{cor}
\noindent

There are $k$ paths of length 2 from a point to itself: going along any edge and then returning along that same edge. This gives $\Delta_2 = k$, and yields Foster's Third:

\begin{cor} \label{}
$$\sum_{w \sim y \sim x \sim v} R_{vw} = 2k^2 + 2k.$$
\end{cor}
\noindent It should be mentioned that in Foster's Third the resistance across all paths of length 3 must be counted, including those which are really paths of length 1 in disguise, i.e. paths of the form $vxvw$ or $vwxw$. However, we can
count the contribution of such paths quite easily. There are $k$ choices for $x$ in each case, however the path $vwvw$ appears in both, so there are a total of $2k-1$ such paths for every $w$ adjacent to $v$.
Thus, the sum over all such paths is simply $(2k-1)\sum_{w \sim v} R_{vw} = 2(2k-1)$, where Foster's First was applied. Thus, Foster's Third can be replaced by the following, which is somewhat simpler to apply:

\begin{cor} \label{}
$$\sum_{\substack{w \sim y \sim x \sim v \\ w \neq x, y \neq v}} R_{vw} = 2k^2 - 2k + 2.$$
\end{cor}
\noindent Let us verify this formula for the square lattice. We do not need to count paths which double back upon themselves, and with that it may be checked that there are three possibilities: we may move three times in the same direction, we may move twice in one direction and once in a perpendicular direction, or we may move three times in three different directions. There are only 4 paths where we move three times in the same direction, corresponding to the 4 possible directions. Moving twice in one direction and once perpendicularly results in a "knight's move", and there are 8 such points to consider; however, each of them may be reached in three different ways, for instance the point $(2,1)$ may be reached by moving twice to the right followed by once up, or once up followed by twice to the right, or once to the right followed by once up followed by once to the right. We see that we have in fact $8 \times 3 =24$ such paths. Finally, if we move three times in three different directions we end up to a point adjacent to the origin. There are 8 ways to do this: 4 choices for the first step, followed by 2 choices for the second step, followed by 1 choice for the final step. Here $k$ is 4, and we verify that

\begin{equation} \label{}
4 \Big(\frac{17}{2}-\frac{24}{\pi}\Big) + 24 \Big(-\frac{1}{2} + \frac{4}{\pi}\Big) + 8\Big(\frac{1}{2}\Big) = 26.
\end{equation}

The formula may also be verified against the values given for the triangular and hexagonal grids, if desired (the hexagonal grid is easy, but beware the triangular one!). Foster's Fourth and higher depends upon the geometry of the graph in question. For instance, if there are no triangles in the graph (such as the square, hexagonal, and truncated square grids discussed earlier), then $\Delta_3 = 0$,
and we have

\begin{cor} \label{} If $G$ has no triangles, then
$$\sum_{w \sim z\sim y \sim x \sim v} R_{vw} = 2k^3 + 2k^2,$$

$$\sum_{\substack{w \sim z\sim y \sim x \sim v\\y \notin \{w,v\}, z \neq x}} R_{vw} = 2k^3 - 4k^2+4k.$$
\end{cor}
\noindent Note that the second equation above is derived from the previous one by subtracting out the 4-paths which are really 2-paths, exactly analogously to what was done earlier in Foster's Third: each 2-path can be realized
as (3k-2) different 4-paths, the resistances across all 2-paths add to $2k$ by Foster's Second, and we obtain a difference of $(3k-2)2k$. As a final example, for the triangular lattice we have $\Delta_3 = 12$, and recall that
$\Delta_2=k=3$ from before. We therefore obtain, on this lattice,

$$\sum_{w \sim z\sim y \sim x \sim v} R_{vw} = 2(6^3) \Big(1 + \frac{6}{6^2} + \frac{12}{6^3}\Big) = 528.$$

\section{Arc transitive graphs and subdivided lattices} \label{split}

Another class of grids upon which Foster's Theorems can be applied is the arc transitive grids, which were defined in the first section. As was noted in Section \ref{pink}, Foster's First and Second easily give the resistances between adjacent points in 1-arc
transitive graphs and points of distance 2 in 2-arc transitive graphs, with the corresponding formulas given there. It is natural, then, to attempt to extend Foster's $r$-th to $r$ larger than 2. However, there are several complications that immediately arise, not least of
which is to give an example of a 3-arc transitive grid. One such example was already discussed in the first section: place a new vertex in the middle of each edge of the hexagonal grid; the resulting grid is no longer vertex
transitive, and in fact is not 2-arc transitive, but is 1-arc and 3-arc transitive.

This observation motivates the consideration of subdivided lattices, which we now define. For a graph $G$, the graph $s(G)$, which we call the {\it subdivision of $G$}, is the graph obtained by placing a new vertex in the middle
of each edge of $G$. It should be mentioned that such lattices are sometimes referred to as {\it decorated lattices}, and have proved useful in various contexts related to statistical mechanics, for instance in counting self-avoiding paths (for example \cite{GL}) and in the analysis of
various ferromagnetic models (for example \cite{fish}). If $v \in s(G)$ and $v$ corresponds to a vertex in $G$ then we will abuse notation somewhat and say $v \in G$, while if $v$ corresponds to an edge of $G$ we will call $v$ a
{\it subdividing vertex}. Unless $G$ is an infinite path, $s(G)$ will not be vertex-transitive,
but if $G$ is edge-transitive and vertex-transitive then $s(G)$ will be edge-transitive, and if $G$ is $2$-arc transitive then $s(G)$ will be $3$-arc transitive. If $G$ is also smallish, then explicit resistances in $s(G)$
up to distance 3 can be calculated by the Foster Theorems, as the next proposition shows. In what follows, $R_{vw}, d(v,w)$ will refer to the resistance and shortest path metrics between points $v,w$ in $G$, while $R'_{vw}, d'(v,w)$
will refer to the corresponding quantities in the subdivision graph $s(G)$.

\begin{prop} \label{emily}
If $G$ is an infinite smallish 1-arc transitive graph of degree $k$, and $v,w$ are adjacent vertices in $s(G)$, then
\begin{equation}
\label{blurb1}
R'_{vw}={{k+2}\over {2k}}.
\end{equation}
Suppose, in addition, that $G$ is 2-arc transitive. If $d'(v,w)=2$ and $v,w \in G$, then
\begin{equation}
\label{blurb3}
R^{'}_{vw}={4 \over k}.
\end{equation}
On the other hand, if $d'(v,w)=2$ and $v,w$ are subdividing vertices, then
\begin{equation}
\label{blurb2}
R^{'}_{vw}={k\over {k-1}}.
\end{equation}
Finally, if $d'(v,w)=3$, then

\begin{equation}
\label{blurb4}
R^{'}_{vw}=\frac{k^2+5k-2}{2k(k-1)}.
\end{equation}
\end{prop}

For the proof, see the appendix. The proposition shows that for the subdivided hexagonal lattice we have $R^{'}_1={{5}\over {6}}$, $R^{'}_{2a}={{4}\over {3}}$, $R^{'}_{2b}={{3}\over {2}}$ and $R^{'}_3={{11}\over {6}}$, where $R^{'}_1$ and $R^{'}_3$ are the resistances between points
of distance 1 and 3 respectively, while $R^{'}_{2a}, R^{'}_{2b}$ are the resistances across the two isomorphism classes of 2-paths. It is worth remarking that the values of the effective resistances found for the subdivided
hexagonal grid are rational, which is not the case for the classic hexagonal grid, where the effective resistances for distances larger than 2 are irrational numbers.

\section{Acknowledgements}

We would like to thank Paul Jung and Tim Garoni for helpful conversations, as well as an anonymous referee for numerous suggestions which have improved the exposition. The first author is grateful for support from Australian Research Council Grants DP0988483 and DE140101201.

\appendix
\section{Appendix}

In this section we give the proofs of Theorem \ref{jess} and Proposition \ref{emily}, and conclude with a few comments on generalizations and the concept of smallishness.

\vski

\noindent {\bf Proof of Theorem \ref{jess}:} We take a large subgraph $G_m$, as before, with the smallish condition on $G$ allowing us to take $m$ to infinity when desired. We apply Foster's $r$-th (with $|V(G_m)| = n$), to obtain

\begin{equation} \label{zhang}
\frac{n}{2} \sum_{j = i_r \sim i_{r-1} \sim \ldots \sim i_0 = v} \frac{R_{vj}}{k^{r-1}} = n \sum_{s=0}^{r-1} P_{vv}^s - r = n \sum_{s=0}^{r-1} \frac{\Delta_s}{k^s} - r.
\end{equation}
We divide both sides by $n$ and let $m \lar \ff$, so that $n$ does as well, and we obtain \rrr{biggie}. \qed

\vski

\noindent {\bf Proof of Proposition \ref{emily}}. As mentioned above, (\ref{blurb1}) is a consequence of Doyle's results, but we can give a quick proof as follows. A large ball in $s(G)$ containing $n$ vertices from $G$ will contain roughly $\frac{nk}{2}$ subdividing vertices
and $nk$ edges. Foster's First shows that the sum of resistances over all edges is $n+\frac{nk}{2}-1$, and thus $R_{vw}$ is approximately $\frac{n+\frac{nk}{2}-1}{nk} \lar {{k+2}\over {2k}}$ as $n \lar \ff$.

Suppose now that $G$ is 2-arc transitive as well. If $d(v,w)=2$ and $v,w \in G$, then $R'_{vw}=2R_{vw}$, because in order to compute $R_{ij}$ in $G^{'}$ when $i$ and $j$ are originally in $G$,
we need to find the voltage at $i$ when a unit current is entered at $i$ and taken out at $j$. But by the way $G^{'}$ is constructed, this voltage is the same if we replace all the pairs of edges between vertices
originally in $G$ with single edges with resistance 2.  Once we have done that, we are back at the original grid $G$ with all edges having resistance 2, with the same unit current flowing from $i$ to $j$,  and thus the voltage
is double that of the original $G$, which was shown to be $\frac{2}{k}$ in \rrr{foster4i}. Now, if $d(v,w)=2$ and $v,w$ are subdividing vertices, then since $G$ is 2-arc transitive $s(G)$ will have exactly two different
isomorphism classes of paths of length 2: those whose middle point is a subdividing vertex, and those whose middle point is a vertex in $G$. In a large ball in $s(G)$ containing $n$ vertices from $G$ there will be
$\frac{nk}{2}$ paths of length 2 whose middle point is a subdividing vertex, and the resistance across each of these paths is approximately $\frac{4}{k}$, as well as $n \frac{k(k-1)}{2}$ paths of length 2 whose middle point is
in $G$, and each of these has a resistance close to $R'_{vw}$. Applying Foster's Second gives

\begin{equation} \label{}
\frac{1}{2}(\frac{nk}{2}) \frac{4}{k} + \frac{1}{k}(n \frac{k(k-1)}{2}) R_{uv} \approx n + \frac{nk}{2} - 2.
\end{equation}
Simplifying and letting $n \lar \ff$ gives \rrr{blurb2}. Finally, if $d(v,w)=3$, then we will apply Foster's Third to a large ball, which says

\begin{equation} \label{marcopolo3}
\sum_{i \in s(G)} \sum_{j \sim u_{2} \sim u_1 \sim i}\frac{R_{ij}}{deg(u_1)deg(u_2)} = 2( tr(P^0) + tr(P^2) - 3).
\end{equation}
Let us calculate the right side first. If our large ball contains $n$ vertices from $G$, then it will also contain roughly $\frac{nk}{2}$ subdividing vertices. Thus, $tr(P^0) \approx n + \frac{nk}{2}$. Each vertex in $G$ is adjacent
to only subdividing vertices in $s(G)$, and vice versa, so if $j \in G$ then $P^2_{jj} = \frac{1}{2}$, while if $j$ is a subdividing vertex then $P^2_{jj} = \frac{1}{k}$. It follows that $tr(P^2) \approx n(\frac{1}{2}) + \frac{nk}{2}(\frac{1}{k}) = n$.
Thus, $2( tr(P^0) + tr(P^2) - 3) \approx 2n + \frac{nk}{2} - 6$. To calculate the left side, we must consider the resistances across all paths, even those of the forms $iu_1ij$ and $iju_2j$. The resistance across all such paths is $R'_1 = \frac{k+2}{2k}$, as was calculated before. If $i$ is in
$G$, then there are $k^2$ paths of the form $iu_1ij$ and $k$ paths of the form $iju_2j$, where $u_2 \neq i$ (but maybe $u_1 = j$). If $i$ is a subdividing vertex, then there are 4 paths of the form $iu_1ij$ and $2(k-1)$ paths of
the form $iju_2j$, where $u_2 \neq i$ (but maybe $u_1 = j$). Also, each genuine 3-arc in $s(G)$ contains a unique point $i$ in $G$ as an endpoint, and there are then $k$ choices for the first step, then only one choice for the second,
and then $k-1$ choices for the third, for a total of $k(k-1)$ paths beginning at $i$. Each of these paths must be counted in the other direction as well, which gives a total of $2nk(k-1)$ three paths in our large ball, with
the resistance across each being approximately $R_{vw}$. Furthermore, in \rrr{marcopolo3}, $deg(u_1)deg(u_2) = 2k$ for every 3-path, since $\{u_1,u_2\}$ must always contain one vertex in $G$ and one subdividing vertex. Combining all
these counts, we have

\begin{equation} \label{}
\frac{1}{2k} \Big( nk(k-1)R_{vw} + \Big(n(k^2 + k) + \frac{nk}{2} (4+ 2(k-1))\Big) R'_1\Big) \approx 2n + \frac{nk}{2} - 6.
\end{equation}
Solving for $R_{vw}$, replacing $R'_1$ with $\frac{k+2}{2k}$, and letting $n \lar \ff$ yields \rrr{blurb4}. \qed

\vski

We now say a few words about Theorem \ref{mpbaby} and the probability theory behind it; a reason for doing this is that a reader interested in extending the methods to smallish resistor networks in which the resistances are not all unity will probably need to understand this. An excellent reference for the interested reader is \cite{DS}. In what follows, $\pi$ will denote the stationary distribution of a random walk $X_t$ upon a finite graph, which in the case of sample random walk (corresponding to unit resistances on all edges) is equal to
$\pi_j = \frac{deg(j)}{2m}$. $T_j$ will denote the first time that $X_t$ hits the vertex $j$, and the standard probabilistic notation $E_iT_j$ will be used for the expectation of $T_j$ under the condition that $X_0 = i$, i.e. the random walk starts at $i$. $P$
will denote the transition matrix of $X_t$ on $G$; that is, $P(X_{t+1} = i| X_t = j) = P_{ij}$, and it follows from this that $P(X_{t+s} = i| X_t = j) = P^s_{ij}$. The
following is the result we use from \cite{P4}.

\begin{thm} \label{fosk}
\begin{equation} \label{shang}
\sum_{i,j \in V(G)} \pi_j P^r_{ji} E_i T_j = \sum_{s=0}^{r-1} tr(P^s) - r.
\end{equation}
\end{thm}
\noindent That this equation can be related to electric resistances is due to a fundamental relation from \cite{Ch}: $E_i T_j + E_j T_i = 2mR_{ij}$. The reversibility of the random walk implies
$\pi_j P^r_{ji} = \pi_i P^r_{ij}$, and using $\pi_j = \frac{deg(j)}{2m}$ we see that both quantities are equal to $\sum_{j \sim v_{r-1} \sim \ldots \sim v_1 \sim i}1/(2mdeg(v_1)deg(v_2) \ldots deg(v_{r-1}))$, where the summation ranges over all paths of length $r$ connecting $i$ to $j$. Using these facts, Theorem \ref{fosk} simplifies to Theorem \ref{mpbaby}. If we are interested in more general cases, the study of grids with variable resistances for instance, then the degree of a vertex $v$ must be defined as $deg(v) = \sum_{y \sim v} C_{vy}$, where $C_{vy} = \frac{1}{R_{vy}}$ is the {\it conductance} between $v$ and $y$. In this case $\pi_j = \frac{deg(j)}{C}$, where $C = \sum_{v \in G} deg(v)$. Equation \rrr{shang} still holds, but applies now to the random walk whose movement on the graph is determined by the resistances along the edges, which is no longer a simple random walk. Furthermore, to pass from \rrr{shang} to resistances, the equation $E_i T_j + E_j T_i = CR_{ij}$ must be used. Results can be obtained in this more general setting, but we have refrained from formulating precise statements due to uncertainty as to applications.

\vski

A final comment on extending Foster's Theorems to infinite graphs. The reader unhappy with the rigor of our method of applying Foster's Theorem to a large ball and then letting the size of the ball go to infinity may also want
to consider the following argument, suggested to us by Rebecca Stones. Smallishness essentially means that our $G$ can be embedded in Euclidean space, and therefore a finite graph which is locally isomorphic to our infinite graph
can be embedded in a compact manifold of the same dimension: for example, for the square lattice, imagine a tessellation of a torus with squares, or of a sphere if one prefers; locally it is identical to the square lattice in $\RR^2$,
Foster's Theorem can be applied directly to it as it is finite, and it is clear that the resistances between points will approach the corresponding ones on the square lattice as we let the size of the squares go to 0. This allows
one to disregard any concerns related to boundary, and shows the validity of the method.

\end{document}